\documentclass{amsart}

\usepackage[utf8]{inputenc} 
\usepackage[T1]{fontenc}    
\usepackage{lmodern}        
\usepackage{amssymb}        
\usepackage{amsthm}         
\usepackage{mathtools}      
\usepackage{mathrsfs}
\usepackage{physics}        
\usepackage{asymptote}      
\usepackage{subcaption}     
\usepackage{xcolor}         
\usepackage[
    hyperfootnotes=false      
]{hyperref}                 
\usepackage{cleveref}       
\usepackage{microtype}      
\usepackage{booktabs}       
\usepackage{tabularx}       
\usepackage[
    maxbibnames=99,           
    maxcitenames=2,           
    backend=biber             
]{biblatex}                 
\usepackage{filecontents}   

\DeclareFieldFormat[article]{pages}{#1}
\DeclareFieldFormat[article]{volume}{\mkbibbold{#1}}
\makeatletter
\def\l@section{\@tocline{1}{0pt}{0pc}{5pc}{}}
\def\l@subsection{\@tocline{2}{0pt}{2.5pc}{5pc}{}}
\makeatother

\definecolor{linkblue}{HTML}{00356B}
\definecolor{linkgold}{HTML}{DA9100}
\definecolor{linkred}{RGB}{159, 29, 53}
\hypersetup{
	colorlinks = true,
	linkcolor = linkred,
	citecolor = linkgold,
	urlcolor = linkblue
}

\theoremstyle{plain}
\newtheorem{theorem}{Theorem}[section]
\crefname{theorem}{Theorem}{Theorems}
\newtheorem{conjecture}[theorem]{Conjecture}
\crefname{conjecture}{Conjecture}{Conjectures}
\newtheorem{proposition}[theorem]{Proposition} 
\crefname{proposition}{Proposition}{Propositions}
\newtheorem{corollary}[theorem]{Corollary} 
\crefname{corollary}{Corollary}{Corollaries}
\newtheorem{lemma}[theorem]{Lemma} 
\crefname{lemma}{Lemma}{Lemmas}
\crefname{ineq}{inequality}{inequalities}

\theoremstyle{definition}

\crefname{example}{Example}{Examples}
\newtheorem{definition}[theorem]{Definition}

\theoremstyle{remark}
\newtheorem*{remark}{Remark}

\crefname{appendix}{Appendix}{Appendices}
\crefname{section}{Section}{Sections}
\crefname{figure}{Figure}{Figures}
\crefname{table}{Table}{Tables}

\newcommand{\field}[1]{\mathbf{#1}}
\newcommand{\R}{\field{R}}

\title{The Optimal Double Bubble for Density $r^p$}

\author[J.\ Hirsch]{Jack Hirsch}
\author[K.\ Li]{Kevin Li}
\author[J.\ Petty]{Jackson Petty}
\author[C.\ Xue]{Christopher Xue}

\thanks{All authors are supported by Yale University through the SUMRY
        fellowship.}

\address{Department of Mathematics, Yale University, New Haven CT, 06510, USA}
\email{jack.hirsch@yale.edu}
\email{k.li@yale.edu}
\email{jackson.petty@yale.edu}
\email{christopher.xue@yale.edu}

\date{\today}

\begin{document}

\begin{asydef}
import geometry;
import patterns;

// Drawing Colors
pen backgroundcolor = white;

picture dot(real Hx=1mm, real Hy=1mm, pen p=currentpen)
{
  picture tiling;
  path g=(0,0)--(Hx,Hy);
  draw(tiling,g,invisible);
  dot(tiling, (0,0), p+linewidth(1));
  return tiling;
}

// Add fill patterns
add("hatch", hatch(1mm, currentpen));
add("dot",dot(1mm, currentpen));

pair bubble(real radius1, real radius2, pen drawingpen = currentpen, int lab = 0) {

	// Draws the appropriate double bubble for bubbles of volume
	// pi * radius1^2 and pi * radius2^2, respectively
	// Calculate the length of the segment connecting the centers of the
	// circles, from the law of cosines. Note that the angle opposite this
	// length formed by the radii is 60 degrees.
	real length = sqrt(radius1^2 + radius2^2 - radius1 * radius2);

	pair origin = (0,0);
	path circle1 = circle(origin, radius1);
	path circle2 = circle((length,0), radius2);

	// Check if the circles have equal curvature. If so, create a straight-line
	// path between their intersection points. Otherwise, draw the bulge as
	// an arc between the intersection points with curvature equal to the 
	// difference of the two circles.
	pair[] i_points = intersectionpoints(circle1, circle2);
	path bulge;
	
	if (radius1 == radius2) {
		// If the radii are equal, the "bulge" will actually be a circle 
		// through infinity, which mathematically is just a line but
		// computationally is an error. In this case, just draw the segment
		// connecting the two points.
		bulge = i_points[0] -- i_points[1];
	} else {
		
		// Calculate the radius of the middle circular arc
		real bulge_radius = 1 / abs(1/radius1 - 1/radius2);

		// Calculate the center of the middle circular arc
		// 	(1) Draw two circles of bulge_radius centered
		//	    at the two singularities. They will intersect
		//	    one another at the two possible centers.
		//	(2) Check which center is correct by comparing
		//	    the original radii.
		path big_circ_1 = circle(i_points[0], bulge_radius);
		path big_circ_2 = circle(i_points[1], bulge_radius);
		pair[] big_i_points = intersectionpoints(big_circ_1, big_circ_2);
		pair big_center;
		if (radius1 > radius2) {
			big_center = big_i_points[1];
		} else {
			big_center = big_i_points[0];
		}
		
		// Calculate the angles from the horizontal between the center of
		// the bulge and the two singularities. Draw an arc at the center
		// using those angles. Since arcsin(...) has a limited range,
		// we must add 180 degrees to theta if the center is on the
		// right, along with reversing the path direction.
		real height = arclength(i_points[0] -- i_points[1]) / 2.0;
		real theta = aSin(height/bulge_radius);
		
		if (radius1 > radius2) {
			bulge = arc(big_center, bulge_radius, -theta + 180, theta + 180);
		} else {
			bulge = arc(big_center, bulge_radius, theta, -theta);
		}
	}

	draw(bulge, drawingpen);

	// Calculate the parameterized time when the og circles intersect
	real[][] i_times = intersections(circle1, circle2);
	real time_1_a = i_times[0][0]; // Time on circle1 of first intersection 
	real time_1_b = i_times[1][0]; // Time on circle1 of second intersection 
	real time_2_a = i_times[0][1]; // Time on circle2 of first intersection 
	real time_2_b = i_times[1][1]; // Time on circle2 of second intersection 

	// The second path is weird to account for the fact that paths always 
	// start at time 0, and we can't run the path in reverse. This then 
	// calculates how long it takes to walk the whole path, and then breaks 
	// the arc up into subpaths.
	// path bubble1 = subpath(circle1, time_1_a, time_1_b) -- bulge -- cycle;
	path bubble1 = subpath(circle1, time_1_a, time_1_b) -- reverse(bulge) -- cycle;
	path bubble2 = subpath(circle2, 0, time_2_a) -- bulge -- subpath(circle2, time_2_b, arctime(circle2, arclength(circle2))) -- cycle;

	filldraw(bubble1, pattern("hatch"), drawingpen);
	filldraw(bubble2, pattern("dot"), drawingpen);

	if (lab == 1) {
		label(Label("$\bm{\Omega_1}$",Fill(backgroundcolor)),origin);
		label(Label("$\bm{\Omega_2}$",Fill(backgroundcolor)),(length, 0));
	}
	
	return i_points[0];
}
\end{asydef}

\begin{abstract}
In 2008 \textcite{Reichardt2007} proved that the optimal Euclidean 
double bubble---the least-perimeter way to enclose and separate two given volumes in $\R^n$---is
three spherical caps meeting along a sphere at 120 degrees. We consider $\R^n$ \emph{with
density} $r^p$, joining the surge of research on manifolds with density after
their appearance in Perelman's 2006 proof of the Poincaré Conjecture.
\textcite{G14} proved that the best single bubble 
is a sphere \emph{through} the origin. We conjecture that the 
best double bubble is the Euclidean solution with 
the singular sphere passing through the origin, for which we have verified equilibrium (first variation or ``first 
derivative'' zero). To prove the exterior of the minimizer connected, it would 
suffice to show that least perimeter is increasing as a function of the 
prescribed areas. We give the first direct proof of such monotonicity in the Euclidean plane. Such arguments were important in the 2002 \emph{Annals} proof \cite{annals} of the double bubble in Euclidean $3$-space.
\end{abstract}

\maketitle
\tableofcontents

\begin{figure}
    \centering
    \includegraphics[width=0.7\textwidth]{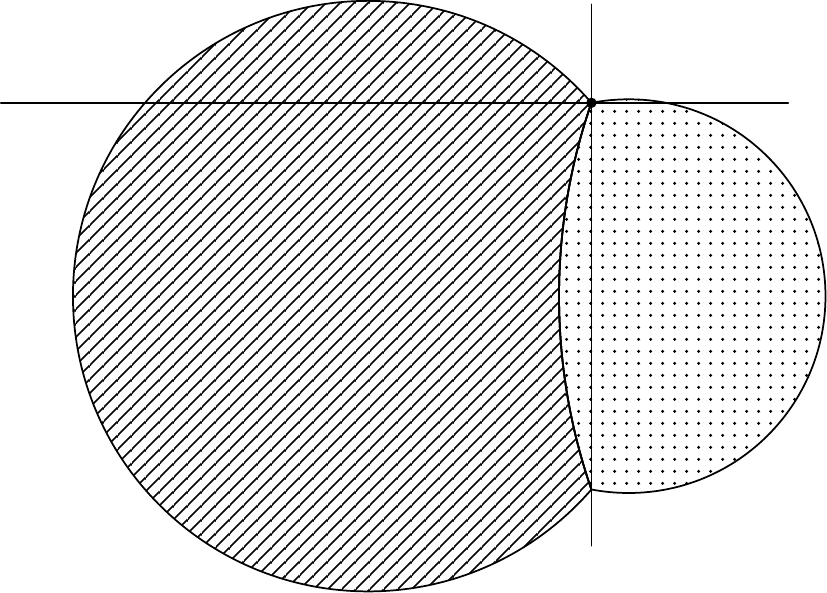}
    \caption{A standard double bubble with vertex at the origin is our 
    conjectured double bubble in the plane with density $r^p$.}
    \label{fig:fdsafdsa}
\end{figure}

\section{Introduction}

The isoperimetric problem is one of the oldest in mathematics. 
It asks for the least-perimeter way to enclose given volume. For a single volume in Euclidean space of any dimension with uniform 
density, the well-known solution is any sphere. In 
Euclidean space with 
density $r^p$, \textcite{G14} found that the solution for a single volume is a sphere
\emph{through} the origin. 
For \emph{two} volumes in
Euclidean space, \textcite{Reichardt2007}
showed that 
the standard double bubble of \cref{fig:double-bubble}, consisting of three spherical caps meeting along a sphere in threes at $120^\circ$ angles, provides an 
isoperimetric cluster. \cref{mainconj} states that the isoperimetric cluster for two volumes in $\R^n$ with density $r^p$ for $p > 0$ is the same Euclidean standard double bubble, with 
the singular sphere passing through the origin, as in 
\cref{fig:fdsafdsa}.

\cref{cor:db-is-equi} verifies equilibrium (first variation or ``first 
derivative'' zero) by scaling arguments and by direct computation. As to 
whether our candidate is the minimizer, it is not even known whether for the 
minimizer each region and the whole cluster are connected. 
Focusing on the 2D case, \cref{prop:nd-ext-connected} notes that to prove the exterior is connected, it 
would suffice to show that the least perimeter 
$P(A_1, A_2)$ for the two areas is increasing in each variable.  
\cref{prop:profile-increasing} gives the first direct proof in the Euclidean 
plane of the ``obvious" but nontrivial fact that $P(A_1, A_2)$ is an increasing
function of the prescribed areas. The original proof of the Euclidean planar double 
bubble by \textcite{foisy1993} finessed the question by considering the 
alternative problem of minimizing perimeter for areas \emph{at least} $A_1$ and
$A_2$, which is obviously nondecreasing. Later \textcite{hutchings} deduced that least perimeter is increasing in higher dimensions from his ingenious proof of 
concavity. Such arguments were important in the 2002 \emph{Annals} proof 
\cite{annals} of the double bubble in Euclidean space.

For our direct proof that $P(A_1, A_2)$ is increasing in the Euclidean plane 
(\cref{prop:profile-increasing}),  we consider the consequences of local minima. In particular, if $P(A_1, A_2)$ is not strictly increasing in $A_1$ for
fixed $A_2$, there is a local minimum never again attained. Because it is a 
local minimum, in a corresponding isoperimetric cluster, the first region has 
zero pressure. Because this minimum is never again attained, the exterior must 
be connected; otherwise a bounded component could be absorbed into the first 
region, increasing $A_1$ and decreasing perimeter. It follows that the dual 
graph has no cycles. Since one can show that components are surrounded by many 
other components as in \cref{fig:branching}, the cluster would have infinitely 
many components, a contradiction of known regularity.

\subsection*{History} Examination of isoperimetric regions in the plane with 
density $r^p$ began in \citeyear{G06} when \textcite{G06} showed that
the isoperimetric solution for a single area in the plane with density $r^p$ is 
a convex set containing the origin. It was something of a surprise when 
\textcite{dahlberg2010} proved that the solution is a circle
through the origin.
In \citeyear{G14} \textcite{G14} extended this result to higher dimensions. In \citeyear{china19} \textcite{china19} studied the $1$-dimensional case,
showing that the best single bubble is an interval with one endpoint at the 
origin and that the best double bubble is a pair of adjacent intervals which 
meet at the origin.
As for the triple bubble, the minimizer in the plane with density $r^p$ cannot 
just be the Euclidean minimizer \cite{wichiramala} with central vertex at the 
origin, because the outer arcs do not have constant generalized curvature.

\subsection*{Acknowledgements}

This paper is a product of the 2019 Yale Summer Undergraduate Research in 
Mathematics (SUMRY) Geometry Group, advised by Frank Morgan. The authors would 
like to thank Morgan for his advice, guidance, and support throughout the 
program, as well as Yale University and the Young Mathematicians Conference for 
providing funding.

\section{Definitions}

\begin{definition}[Density Function]
Given a smooth Riemannian manifold $M$, a \emph{density} on $M$ is a positive 
continuous function, which weights each point $p$ in $M$ with a certain mass 
$f(p)$. Given a region $\Omega \subset M$ with piecewise smooth boundary, the weighted volume and 
perimeter of $\Omega$ are given by
\[
V(\Omega) = \int_{\Omega} f \dd{V_0} \quad\text{and}\quad A(\Omega) = 
\int_{\partial \Omega} f \dd{P_0},
\]
where $\dd{V_0}$ and $\dd{P_0}$ denote Euclidean volume and perimeter. We may also refer to the perimeter of $\Omega$ as the perimeter of its boundary.
\end{definition}

\begin{definition}[Isoperimetric Region]
A region $\Omega \subset M$ is 
\emph{isoperimetric} if it has the smallest weighted perimeter of all regions 
with the same weighted volume. The boundary of an isoperimetric region is also 
called isoperimetric.
\end{definition}

We can generalize the idea of an isoperimetric region by considering two or more
volumes.
\begin{definition}[Isoperimetric Cluster]
An isoperimetric cluster is a set of disjoint open regions~$\Omega_i$ of 
volume~$V_i$ such that 
the perimeter of the union of the boundaries is minimized.
\end{definition}

To provide an example of the concepts we have introduced, consider the 
isoperimetric solution for a single unit volume in $\R^n$ with constant density 
$1$. The solution is simply a sphere.

\begin{figure}
    \centering
    \begin{subfigure}[c]{0.3\textwidth}
        \centering
        \includegraphics[height=1.3in]{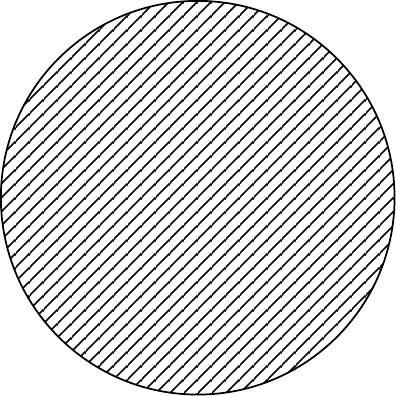}
        \caption{$f = 1$}
        \label{fig:unit-circle}
    \end{subfigure}%
    \hspace{2cm}
    \begin{subfigure}[c]{0.3\textwidth}
        \centering
        \includegraphics[height=1.457in]{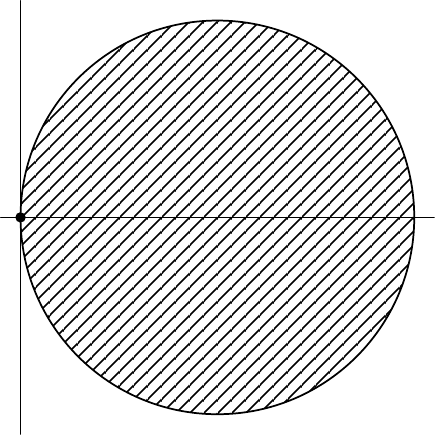}
        \caption{$f = r^2$}
        \label{fig:unit-circle-r2}
    \end{subfigure}
    \caption{Known single-volume isoperimetric solutions. In the Euclidean plane,
             it is any circle of the prescribed area; for density $r^p$, it is a
             circle \emph{through the origin}.}
    \label{fig:known-solution}
\end{figure}

For density $r^p$, the solution in the plane is a circle passing through the origin 
\cite[Thm.\ 3.16]{dahlberg2010}, as shown in \cref{fig:known-solution}; in higher dimensions, the solution is a sphere passing through the origin \cite[Thm.\ 3.3]{G14}.

\textcite{annals} proved in \citeyear{annals} that the isoperimetric solution for two volumes in Euclidean space with constant density is the standard double bubble, so called because of how soap bubbles combine in three-dimensional space, as in \cref{fig:double-bubble}. The standard double bubble illustrates the existence, boundedness, and regularity theorems:

\begin{figure}
    \begin{subfigure}[t]{0.3\textwidth}
        \centering
        \includegraphics[height=1.3in]{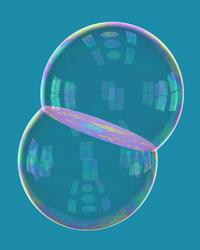}
    \end{subfigure}
    \hspace{2cm}
    \begin{subfigure}[t]{0.3\textwidth}
        \centering
        \includegraphics[height=1.3in]{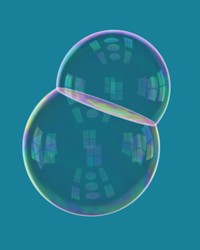}
    \end{subfigure}
    \caption{The standard double bubbles for volumes $V_1 = V_2$ and $V_1 > 
             V_2$. John M. Sullivan, \url{http://torus.math.uiuc.edu/jms/Images/double/}, used with permission}
    \label{fig:double-bubble}
\end{figure}

\begin{lemma}\label{lem:escaping-area}
Consider $\R^n$ with radial non-decreasing density $f$ such that 
$f(r) \to \infty$ as $r \to \infty$. If a sequence of clusters  $\Omega_i$ with uniformly 
bounded perimeter converge to $\Omega$, there is no loss of volume at
infinity in the limit. 
\end{lemma}

\begin{proof}
This fact is shown almost identically in the proof of one region by \textcite[Thm.\ 2.1]{rosales}. Consider each region of the sequence of clusters separately to obtain a sequence of regions of fixed volume and uniformly bounded perimeters. The proof of \textcite[Thm.\ 2.1]{rosales} implies that there is no loss in volume at infinity for each region and therefore for the whole cluster.
\end{proof}

\begin{theorem}[Existence]\label{thm:existence}
Consider $\R^n$ endowed with a nondecreasing radial density $f$ such that 
$f(r) \to \infty$ as $r \to \infty$. Given volumes $V_1, \dots, V_n$, there 
exists an isoperimetric cluster that separates and encloses the given volumes. 
\end{theorem}

\begin{proof}
The proof is almost identical to the proof of one region by \textcite[Thm.\ 
2.5]{rosales}, because their argument does not depend on the number of regions. 
We give the generalization of the proof here to multiple regions for the 
convenience of the reader. Consider a sequence of clusters enclosing and 
separating volumes $V_1, \dots, V_n$ such that their perimeter  approaches 
$I(V_1, \dots, V_n)$ and is less than $I(V_1, \dotsc, V_n) + 1$. By the 
Compactness Theorem \cite[Sect.\ 9.1]{morgan}, we may assume this sequence 
converges. By \cref{lem:escaping-area}, there is no loss of volume at infinity, so
the limit gives the isoperimetric region. 
\end{proof}

\begin{proposition}[Boundedness] \label{prop:boundedness}
In $\R^n$ with radial non-decreasing $C^1$ density $f$, a perimeter-minimizing cluster is bounded.
\end{proposition}

\begin{proof}
The proof follows \textcite[Thm.\ 5.9]{Morgan2013}, generalizing their proof to clusters. Some analogous details of their proof are omitted for brevity. Note that while their proof is emphasized for $n \ge 3,$ it nonetheless applies for $n=2$. Assume to the contrary that $E$ is an unbounded isoperimetric set with unbounded regions $E_1,\ldots,E_b$ and bounded regions $E_{b+1},\ldots,E_n.$ Each region has bounded perimeter.

Let $B(r)$ and $S(r)$ be the closed ball and sphere of radius $r$, and let $\mathscr{H}_f^{k}$ denote the $k$-dimensional Hausdorff measure in $\R^n$ with density $f$. Define
\begin{align*}
    E_i(r)&:=E_i \cap B(r) &E_i^r := E_i \cap S(r) \\
P_i(r) &:= \mathscr{H}_f^{n-1}(\partial E_i\backslash B(r)) &V_i(r) := \mathscr{H}_f^{n}(E_i\backslash B(r)).
\end{align*}
 Since $E_1$ is unbounded, the proof of \cite[Thm.\ 4.3]{Morgan2013} shows
\[P(E_1(r)) < P(E_1)\]
which, as $P(E_1(r)) = P(E_1) - P_1(r) + \mathscr{H}^{n-1}_f(E_1^r),$ is equivalent to
\[P_1(r) > \mathscr{H}^{n-1}_f(E_1^r).\]


After a careful application of the isoperimetric inequality on the $(n-1)$-sphere and some manipulations, the above inequality yields for sufficiently large $r$
\begin{equation}
    \label[ineq]{iso}
    P_1(r)^{\frac{n}{n-1}} \ge cV_1(r)
\end{equation}
for some positive constant $c$.

Note that there must be some unbounded region with a component $C$ which borders the exterior. Without loss of generality, let this be the unbounded region $E_1.$ Pick an $R$ such that $B(R)$ completely contains the bounded components of $E$ and such that \cref{iso} holds true for any $r>R$, and such that the part of $C$ which borders the exterior has nonzero measure $\mathscr{H}_f^{n-1}$ in $B(R)$. There exists a constant $\overline{\varepsilon}$ such that, for any $0 < \varepsilon < \overline{\varepsilon},$ it is possible to make small variations to  $C$ along the exterior and replace the set $E_1$ with another set $E^{\varepsilon}$ such that
\begin{align*}E^{\varepsilon}\backslash B(R) = E_1\backslash B(R), \text{ }V(E^{\varepsilon})=V(E_1)+\varepsilon, \text{ } P(E^{\varepsilon}) \le P(E_1) + \varepsilon(H(E_1)+1),\end{align*}
where $H$ is the generalized mean curvature (\ref{def:gen_curve}), well defined because the density $f$ is $C^1$.

Now, for any $r$ sufficiently large, set $\varepsilon = V_1(r) < \overline{\varepsilon}$ and $\widetilde{E} = E^{\varepsilon} \cap B(r)$. By construction, \[V(\widetilde{E}) = V(E_1)\] and
\begin{align*}P(\widetilde{E}) &= P(E^{\varepsilon}) - P_1(r)+\mathscr{H}_f^{n-1}(E_1^r)\\& \le P(E_1)+\varepsilon(H(E_1)+1) - P_1(r) + \mathscr{H}_f^{n-1}(E_1^r).\end{align*}

Since the volume of $\widetilde{E}$ equals the corresponding volume $V(E_1)$ and the variation was made only to $E_1$ and only along the exterior, i.e., no shared perimeter changes, it must be that
\[ P(E_1) \le  P(\widetilde{E})\]
as $E$ is isoperimetric for its given volumes.



Hence, taking $\varepsilon$ arbitrarily small by picking arbitrarily large $r >R$ and using \cref{iso}, it follows that
\[\mathscr{H}_f^{n-1}(E_1^r) \ge c^{\frac{n-1}{n}}V_1(r)^{\frac{n-1}{n}},\]
which is equivalent to
\[-\frac{\partial}{\partial r}\left(V_1(r)^\frac{1}{n}\right) \ge \frac{c^{\frac{n-1}{n}}}{n},\]
contradicting the fact that $V_1(r)>0$ for all $r$.

\end{proof}
\begin{definition}[Generalized Curvature]\label{def:gen_curve}
In $\R^2$ with density $f$, the 
generalized curvature $\kappa_f$ of a curve with inward-pointing unit normal $N$ is given by the formula
\[
\kappa_f = \kappa_0 - \pdv{\log f}{N},
\]
where $\kappa_0$ is the (unweighted) geodesic curvature. This comes from the first variation formula, so that generalized curvature has the interpretation as minus the perimeter cost $\dd{P}\!/\!\dd{A}$ of moving area across the curve, and constant generalized mean curvature is the equilibrium condition $\dd{P} = 0$ if $\dd{A} = 0$ (see \textcite[Sect.\ 3]{rosales}).

More generally, for a smooth open region $\Omega$ in $\R^{n + 1}$ with boundary $\Sigma$ with smooth density $f = e^\psi$, we can define the generalized mean curvature to be 
\[H_f = H_0 - \langle \nabla \psi, N \rangle,\]
where $N$ is the inward normal unit vector to $\Sigma$, and $H_0$ is the Euclidean mean curvature (sum of principal curvatures) with respect to $N$. 
\end{definition}

\begin{theorem}[Regularity]\label{thm:regularity}
An isoperimetric cluster in $\R^2$ with smooth density consists of smooth 
constant-generalized-curvature curves meeting in threes at $120^\circ$. The sum 
of the curvatures encountered along a generic closed path is $0$.


\end{theorem}
\begin{proof}An isoperimetric cluster is a so-called 
(\textbf{M},$Cr^\alpha,\delta$)-minimal set, and therefore consists of curves 
meeting in threes at $120^\circ$ (see \textcite[Sect.\ 13.10]{morgan}). The rest is 
the equilibrium conditions (see \textcite[Sect.\ 3]{rosales}).
\end{proof}

For regularity in higher dimensions, see \textcite[Sect.\ 13.10]{morgan} for a detailed discussion.

\begin{remark}
    Consider, in the plane of density $r^p,$ a circle $C$ of radius $R$ centered at $(x_0, y_0)$. At some point $(a, b) \in C,$ the normal vector is $\frac{1}{R}(a-x_0, b-y_0).$ If $(a,b) \neq (0,0),$ the generalized curvature is
    \begin{align*}
        \kappa_{r^p} &= \kappa_0 - \pdv{\log f}{n} = \frac{1}{R} - \frac{p}{2}\pdv{\log (a^2+b^2)}{n} \\[1ex] 
        &= \frac{1}{R} - \frac{p}{R}\frac{a(a-x_0) + b(b-y_0)}{a^2+b^2} \\[1ex]
        &= \frac{1}{R} - \frac{p}{R}\left(1 - \frac{ax_0+by_0}{a^2+b^2}\right).
    \end{align*}
    From \[(a-x_0)^2 + (b-y_0)^2 = R^2\] we find \[a^2 + b^2 + x_0^2 + y_0^2 = R^2+2(ax_0 + by_0),\] therefore $(ax_0+by_0)/(a^2+b^2)$ is constant if and only if $(x_0^2+y_0^2 - R^2)/(a^2+b^2)$ is. Evidently, this happens if and only if either $x_0^2+y_0^2 = R^2$ or $a^2+b^2$ is constant. In other words, $C$ has constant generalized curvature if and only if it either passes through or is centered at the origin.
    
    This result extends straightforwardly to $\R^n$: the spheres in $\R^n$ with density $r^p$ with constant generalized curvature are precisely those passing through or centered at the origin.
\end{remark}
\section{Double Bubble in density \texorpdfstring{$r^p$}{rp}}


We conjecture that the isoperimetric cluster for two regions in $\R^n$ with density $r^p$ has the exact same shape as the Euclidean standard double bubble, 
but with the singular sphere passing through the origin. Notice that every cap is now part of a sphere through the origin, proved by \textcite{G14} to be the best 
\emph{single} bubble. 

\begin{proposition}
For any two given volumes, there is a standard double bubble with singular sphere passing through the origin, unique up to rotation.
\end{proposition}

\begin{proof}
This proof follows directly from \textcite[Prop.\ 14.1]{morgan}, the existence of unique standard bubbles in unit density.
\end{proof}
\begin{conjecture}
\label{mainconj}
Consider $\R^n$ with density $r^p$ for positive $p$. The isoperimetric solution for two regions in space is the standard double bubble with singular sphere passing through the origin, unique up to rotation.
\end{conjecture}

The proof (\cref{cor:db-is-equi}) that our candidate is in equilibrium (first 
variation zero) will require the following scaling lemma.

\begin{lemma}\label{lem:curvescale}
In the space $\R^n$ with density $r^p$, if a surface is scaled by $\lambda$ about the 
origin, then the generalized curvature is scaled by $1/\lambda$.
\end{lemma}
\begin{proof}
In space with density $r^p$, perimeter is scaled by $\lambda^{p+n-1}$, and 
volume is scaled by $\lambda^{p+n}$. Since generalized curvature has the 
interpretation of $\dd{P}\!/\!\dd{V}$, it is scaled by $1/\lambda$. 
\end{proof}

\begin{corollary} \label{cor:db-is-equi}
The standard double bubble in $\R^n$ with density $r^p$ for some $p > 0$ and 
singular sphere passing through the origin is in equilibrium. 
\end{corollary}

\begin{proof}
\textcite{Reichardt2007} showed that the standard double bubble in $\R^n$ with unit density is isoperimetric, in particular in equilibrium. Thus, the three spherical caps meet at 120 degrees, have constant Euclidean curvature, and the sum of the Euclidean curvatures encountered along a generic closed path is 0. Observe that the spherical caps also have constant generalized curvature since they all pass through the origin. By \cref{lem:curvescale}, their generalized curvatures are in proportion to their inverse radii, i.e., Euclidean curvature. It follows that the sum of the generalized curvatures encountered along a generic closed path must be 0.
\end{proof}

\begin{proposition}\label{prop:pospressure}
In a bounded isoperimetric cluster in $\R^n$ with non-decreasing radial 
density, the region farthest from the origin must have positive pressure. 
\end{proposition}
\begin{proof}
Since at the point farthest from the origin the cluster lies in a halfspace, the tangent cone must be a hyperplane and the cluster must be regular by \textcite{Morgan2003RegularityOI}, with normal vector pointing
towards the origin. Also, at the point farthest from the origin, the unweighted Euclidean curvature is positive. Since the log of the density is radially non-decreasing, the generalized 
curvature also must be positive, and hence the region must have positive pressure. 
\end{proof}

\section{Geodesics in plane with density 
\texorpdfstring{$r^p$}{rp}}\label{sec:geodesics-in-plane}

Geodesics in the plane with density $r^p$ can be completely analyzed by mapping 
the plane with density to a Euclidean cone with area density.

\begin{proposition}\label{areadensitymap}
The conformal map $w= z^{p+1}/(p+1)$ takes the plane with area and perimeter 
density $r^p$ to a Euclidean cone with angle $(p+1)\pi$ about the origin, with 
area density $r^{-p}\sim|w|^{-p/(p+1)}$
(and perimeter density $1$). 
\end{proposition}
\begin{proof}
Since the derivative $z^p$ has modulus $r^p$, the image perimeter density is 
$r^{-2p}r^p=r^{-p}$ and the image area density is $r^{-p} r^p=1$.  
\end{proof}
Note that in the image, a geodesic is either a straight line or two straight 
lines meeting at the origin.

\begin{corollary}
In the plane with density $r^p$, the unique geodesic from any point to the 
origin is the straight line. For two points with $\Delta \theta$ at least 
$\pi/2(p+1)$, the unique geodesic consists of two lines to the origin. For two 
points with $\Delta \theta$ less than $\pi/2(p+1)$, there is a unique geodesic 
corresponding to a straight line segment in the Euclidean cone.
\end{corollary}

\section{Properties of the Isoperimetric Function}

Understanding the isoperimetric function, how least perimeter depends on volume, 
has important consequences for the shape of minimizers.
We begin with some preliminary results about scaling.
As noted by \textcite[Sect.\ 3.6]{dahlberg2010} for the plane 
and mentioned in our proof of \cref{lem:curvescale}, the density $r^p$ has nice scaling 
properties. If a cluster $\Omega$ has perimeter $P$ and volume $V$, then 
$\lambda\Omega$ has perimeter $\lambda^{p+n-1}P$ and volume $\lambda^{p+n}V$.

\begin{lemma}
In Euclidean space with density $r^p$, if a cluster $\Omega$ is scaled such that the 
volume is scaled by $\lambda$, then the perimeter is scaled by 
$\lambda^{\frac{p + n - 1}{p + n}}$.
\label{lem:scalearea}
\end{lemma}

\begin{proof}
When the cluster $\Omega$ is scaled to $\lambda^{\frac{1}{p + n}} \Omega$, the 
volume is scaled by $\lambda$ and the perimeter is scaled by 
$\lambda^{\frac{p + n-1}{p + n}}$.
\end{proof}

\begin{definition}[Isoperimetric Function]
The isoperimetric function $I(V_1, V_2)$ as the least perimeter to enclose and 
separate volumes $V_1$ and $V_2$.
In our applications, minimizers exist. In general, $I$ would be defined as an 
infimum. 
\end{definition}

The isoperimetric function $I$ reflects the nice scaling properties of density 
$r^p$.

\begin{proposition}
In $\R^n$ with density $r^p$, for any volumes $v$ and $w$, 
\[
I(\lambda v, \lambda w) = \lambda^{\frac{p + n - 1}{p + n}}I(v, w).
\]
\end{proposition}
\begin{proof}
By \cref{lem:scalearea}, 
$I(\lambda v, \lambda w) \le \lambda^{\frac{p +n- 1}{p + n}} I(v, w)$. Reapplying 
the lemma with $1/\lambda$ yields the opposite inequality.
\end{proof}

The following proposition that the isoperimetric profile is 
continuous is by no means clear for spaces of infinite measure. Indeed, 
\textcite{nardulli} and \textcite{papasoglu} give examples of (noncompact) two- 
and three-dimensional Riemannian manifolds with discontinuous isoperimetric 
profile.

\begin{proposition}
In $\R^n$ with density $r^p$, the isoperimetric profile $I(v,w)$ is 
continuous. 
\end{proposition}
\begin{proof}
To prove upper semicontinuity, just note that small changes in volume can be 
attained by a small change in perimeter. For lower semicontinuity, consider a 
sequence of volumes $(v_i, w_i) \to (v, w)$. Let $\Omega_i$ be a isoperimetric 
cluster of volume $(v_i, w_i)$. By the Compactness Theorem 
\cite[Sect.\ 9.1]{morgan}, we may assume that $\Omega_{i} \to \Omega$, and by 
\cref{lem:escaping-area}, volume does not escape to infinity, so $\Omega$ encloses
and separates volumes $v$ and $w$. By the lower semicontinuity property 
\cite[Ex.\ 4.22]{morgan} $A(\Omega_{i}) \le \liminf_{i \to \infty} I(v_i, w_i)$.
Since $I(v, w)$ is the perimeter of the isoperimetric cluster, we must have 
$I(v, w) \le \liminf_{i \to \infty} I(v_i, w_i)$. Therefore $I$ is lower 
semi\-continuous, and hence continuous.
\end{proof}

Properties of the isoperimetric profile imply connectivity properties of an 
isoperimetric cluster. \cref{prop:subadditive} gives the trivial implication 
that $I$ subadditive implies that the cluster is connected.

\begin{proposition}\label{prop:subadditive}
Consider a Riemannian manifold with density where an isoperimetric cluster 
exists for all volumes. If the isoperimetric profile is strictly subadditive, 
then any isoperimetric cluster is connected.
\end{proposition}
\begin{proof}
Suppose the isoperimetric cluster of volumes $(v, w)$ is not connected. Then the
cluster can be separated into two disjoint clusters, one with volumes 
$(v_1, w_1)$ and another with volumes $(v_2, w_2)$ with $v_1 + v_2 = v$ and 
$w_1 + w_2 = w$. Then 
\[
I(v, w) \le I(v_1, w_1) + I(v_2, w_2),
\]
which contradicts strict subadditivity. 
\end{proof}

The following proposition proves in a general context that $I$ increasing 
implies that the exterior is connected.

\begin{proposition} \label{prop:nd-ext-connected}
Consider $\R^n$ with radial density $f(r)$ such that $\liminf_{r \to \infty} f(r) > 0$.
If the isoperimetric profile is non-decreasing, the exterior of an isoperimetric
cluster is connected.
\end{proposition}

\begin{proof}
Since $\liminf f(r)$ does not vanish, an unbounded hypersurface yields infinite 
weighted perimeter. Therefore the unbounded component of the exterior is 
connected. If there is a bounded component, simply absorb it into an adjacent 
region, which decreases the perimeter and increases volume, a contradiction of the assumption that 
the isoperimetric profile is non-decreasing. 
\end{proof}

On the other hand, focusing on the plane, the following proposition shows that if the isoperimetric profile is not increasing, at least one of the regions is not connected.

\begin{proposition}
In the plane with density $r^p$, if the isoperimetric profile $I$ is not 
(strictly) increasing in each variable, then there exists an isoperimetric 
cluster such that the region farthest from the origin has at least two 
components.
\end{proposition}

\begin{proof}
Since $I$ is not increasing, there exist $v_0$ and $w_0$ such that say 
$I(v_0, w_0)$ is a local minimum of $I_{v_0}(w) = I(v_0, w)$ and
\begin{equation}
    I(v_0, w_0) \le I(v_0, w) \qquad \text{for all $w > w_0$.} 
    \tag{\text{$\star$}} \label{eq:profile-increasing}
\end{equation}
Since $w_0$ is a local minimum, the second region $R_2$ must have $0$ pressure. 
The image of each component of $R_2$ under the map of \cref{areadensitymap} to 
the flat cone with only area density is bounded by geodesics and negative 
curvature curves meeting at $120^\circ$, bounding alternately $R_1$ and the 
exterior. If it does not pass through the origin, it has at least eight edges. 
Since regularity does not hold at the origin, where the density is $0$, a 
geodesic could turn at a small angle there, but it still has at least four 
edges, two of  which border $R_1$. To see that they are different components of 
$R_1$, note that the exterior cannot have a bounded component, because because 
such a component could be absorbed into $R_2$, contradicting 
\eqref{eq:profile-increasing}. Therefore the component is bounded by at least 
two distinct components of $R_1$. 
\end{proof}

The following proposition proves for the Euclidean plane the ``obvious'' but nontrivial fact that least 
perimeter $I(v,w)$ is an increasing function of the prescribed areas. The 
original proof of the Euclidean double bubble by \textcite{foisy1993} finessed 
the question by considering the alternative problem of minimizing perimeter for 
areas \emph{at least} $v$ and $w$, which is obviously nondecreasing. Later 
\textcite{hutchings} deduced $I$ increasing in higher dimensions from his 
ingenious proof of $I$ concave.

\begin{proposition} \label{prop:profile-increasing}
In the Euclidean plane, the isoperimetric profile $I(v, w)$ is (strictly) 
increasing in each variable. 
\end{proposition}

\begin{figure}
\centering
\includegraphics[width=0.5\textwidth]{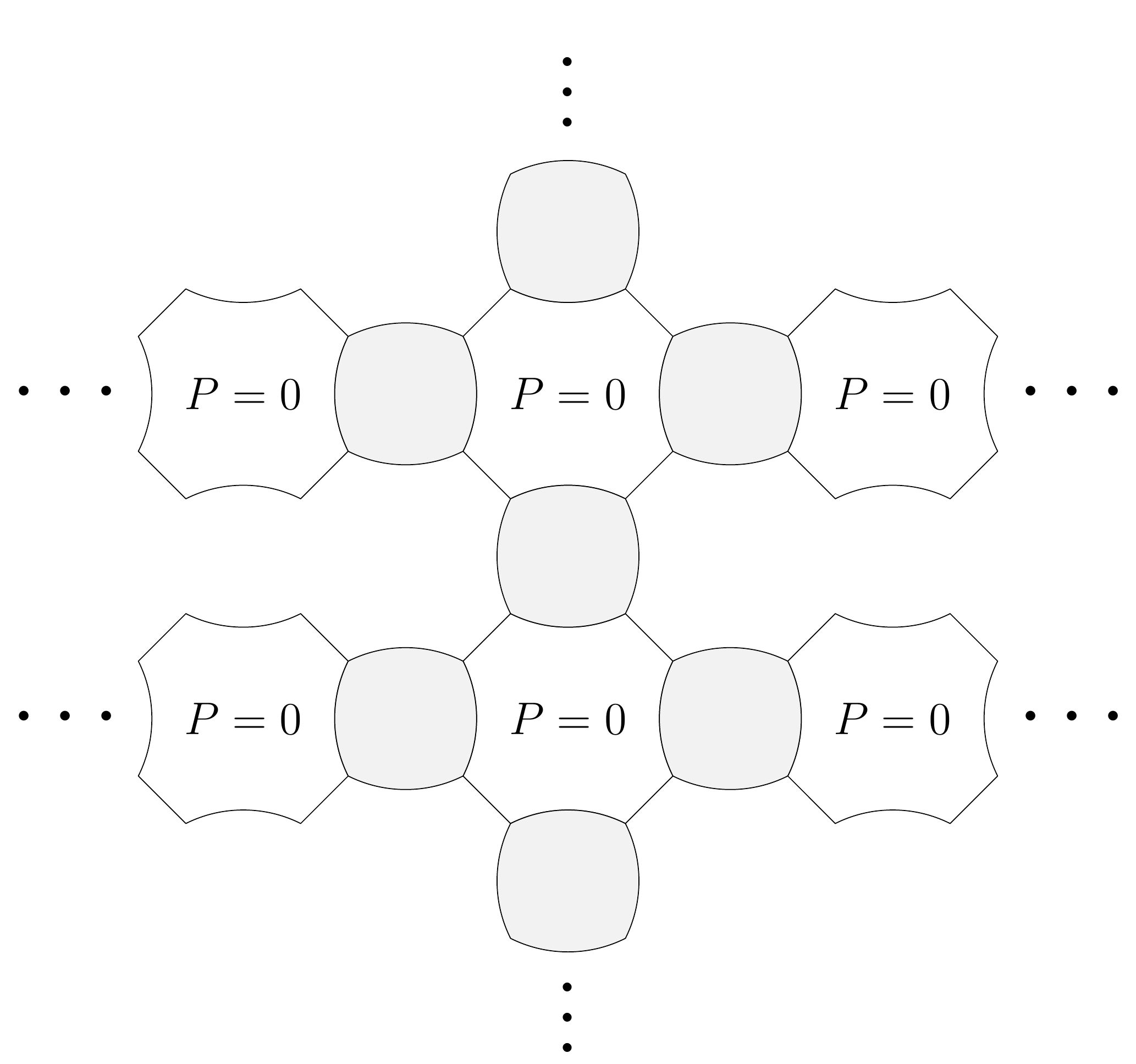}
\caption{Our direct proof that the isoperimetric function $I(A_1, A_2)$ in the 
Euclidean plane is increasing shows that otherwise there would be a region of 
pressure $0$ and infinite branching.}
\label{fig:branching}
\end{figure}

\begin{figure}
    \includegraphics{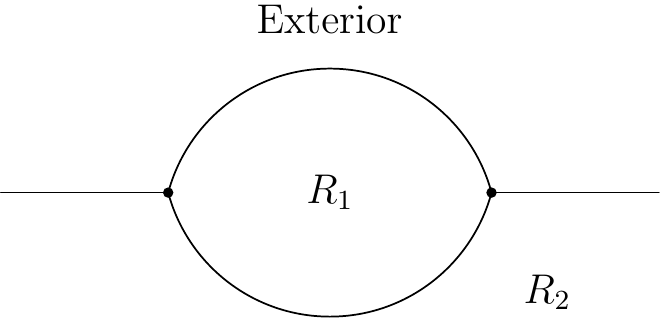}
    \caption{If a component of $R_1$ has just two edges, then the adjacent two  
             edges bounding $R_2$ are two collinear geodesics.}
    \label{fig:medial-bubble}
\end{figure}

\begin{proof}
If not, there exists a $v_0$ such that $I_{v_0}(w) = I(v_0, w)$ is not 
increasing. Since $I_{v_0}(w) \to \infty$ as $w \to \infty$ and is continuous, 
there exists a $w_0$ such that $I_{v_0}(w_0)$ is a local minimum and 
$I_{v_0}(w_0) < I_{v_0}(w)$ for all $w > w_0$. Let $\Omega$ be an isoperimetric 
cluster of areas $v_0$ and $w_0$, and let $R_1$ and $R_2$ denote the regions of 
areas $v_0$ and $w_0$ respectively. The second region $R_2$ must have zero 
pressure, since otherwise, it is possible to decrease perimeter while changing 
area. On the other hand, by \cref{prop:pospressure} the region farthest from the
origin must have positive pressure and hence must be $R_1$. 

The exterior of $\Omega$ is connected; otherwise, since the cluster is bounded 
(\cref{prop:boundedness}), there would be a bounded component of the exterior, 
which could be absorbed into $R_2$, contradicting $I_{v_0}(w_0) < I_{v_0}(w)$ 
for all $w > w_0$. Therefore the dual graph of $\Omega$, with a labeled vertex 
for each component of $R_1$ and of $R_2$, does not have any cycles. Since as in 
\cref{fig:branching} a component of $R_2$ is bounded by alternating geodesic and
strictly concave segments meeting at $120^\circ$, it has at least eight edges. 
Since a component of $R_1$ is convex with $120^\circ$ angles, it must have two 
or four edges (alternately shared with $R_2$ and the exterior). If it has two 
edges as in \cref{fig:medial-bubble}, the two adjacent geodesics are collinear. 
Hence at least two components of $R_1$ on the boundary of every component of 
$R_2$ both have four edges. Since the dual graph has no cycles, starting at a 
component of $R_2$, moving to an adjacent component of $R_1$ with four edges, 
moving to the other adjacent component of $R_2$, moving to another adjacent 
component of $R_1$, etc., would yield infinitely many components, a 
contradiction of boundedness and regularity (\cref{prop:boundedness}, 
\cref{thm:regularity}).
\end{proof}

\begin{remark}
Note that this argument does not extend to $\R^3$, since the dual graph may 
contain cycles even if the exterior is connected. 
\end{remark}

\begin{remark}
For an isoperimetric double bubble in $\R^3$ with unit density, 
\textcite[Sect.\ 4.5]{hutchings} proved the much stronger result that there are 
at most three components. In particular, the larger region has exactly one 
component, and the smaller region has at most two components. 
\end{remark}

\section{Connectedness}

Although \emph{a priori} we do not know that an isoperimetric cluster is 
connected, we can say something about what multiple components would look like. 

\begin{lemma}\label{lemma:annulus}
In $\R^n$ with radial density, components of an isoperimetric cluster must 
lie in disjoint open hyperspherical shells. 
\end{lemma}

\begin{proof}
If components are not in disjoint open hyperspherical shells, then one can rotate a component 
about the origin until it contacts another component, contradicting regularity 
(\cref{thm:regularity}). 
\end{proof}

In general, a conformal map takes a surface with density to a surface with 
different area and perimeter densities. For the right conformal map, however, 
one of the densities could be made to be $1$. With unit area density it is 
easier to find transformations that preserve area.

Now we focus on the two-dimensional problem. A conformal map takes a surface with density to a surface with 
different area and perimeter densities. For the right conformal map, however, 
one of the densities can be made to be $1$. With unit area density it is 
easier to find transformations that preserve area.

\begin{proposition}\label{prop:perimdensitymap}
The conformal map 
\[
w = \frac{2}{p + 2}z^{\frac{p + 2}{2}}
\]
takes the plane with area and perimeter density $r^p$ to a Euclidean cone with 
perimeter density $r^{p/2}\sim|w|^{p/(p+2)}$ (and area density $1$).

\end{proposition}
\begin{proof}
Since the derivative $z^{p/2}$ has modulus $r^{p/2}$, the image perimeter 
density is $r^{-p/2}r^p=r^{p/2}\sim|w|^{p/(p+2)}$ and the image area 
density is $r^{-p} r^p=1$.  
\end{proof}

The following lemma gives a nice map that preserves area.

\begin{lemma}\label{lemma:areapreserved}
Consider an open set $U$ in the plane with area density 1 such that $U$ is 
outside some ball $B(0, \sqrt{\epsilon})$. Let $\varphi_\epsilon(r) := \sqrt{r^2-\epsilon}.$ The polar map
\[
\Phi_\epsilon\colon (r, \theta) \mapsto (\varphi_\epsilon(r), \theta)
\] 
preserves the area of $U$. 
\end{lemma}

\begin{proof}
Note that $\Phi_\epsilon\colon U \to \R^2$ is injective. A computation shows that the 
determinant of the Jacobian $\det(D\Phi_\epsilon) = 1$, so the area of $U$ is 
preserved. 
\end{proof}

\begin{remark}
Given $\epsilon > 0$, the maps $\Phi_\alpha$ for $\alpha \le \epsilon$ are actually the only radially symmetric differentiable maps that preserve 
area outside $B(0, \sqrt{\epsilon}).$ Indeed, suppose $\Phi(r, \theta) = (\varphi(r), \theta)$ preserves area outside that ball. For every open set $V \subseteq \mathbb{R}^2 \setminus B(0, \sqrt{\epsilon})$,
\[
    \int_V r \dd{r} \dd{\theta} = \int_V \varphi(r) \det(D\varphi) \dd{r} \dd{\theta} 
                         = \int_V \varphi(r) \varphi'(r) \dd{r} \dd{\theta}.
\]
Thus, $\varphi$ must satisfy 
$$r = \varphi(r) \varphi'(r) =  \frac{1}{2}(\varphi(r)^2)'$$ for almost all $r \ge \sqrt{\epsilon}.$ 
One can extend this equation to all $r \ge \sqrt{\epsilon}$ by continuity and solve to conclude that $\varphi$ takes the form $\varphi(r) =\sqrt{r^2 - \alpha}$ for some $\alpha \le \epsilon.$
\end{remark}

Next we consider how the map $\varphi_\epsilon$ affects the perimeter. 

\begin{lemma}\label{lemma:lessperim}
Consider a smooth curve in the plane with perimeter density $r^{k}$ with 
$k > 1$, outside some ball $B(0,\sqrt{\epsilon})$. The map $\Phi_\epsilon$ strictly decreases the the length of the curve. 
\end{lemma}
\begin{proof}
Note that $\Phi_\epsilon$ clearly decreases the length of infinitesimal tangential 
elements. Therefore it suffices to consider an infinitesimal radial element at 
$r$. The Euclidean length is scaled by 
\[
\lambda = \dv{r} \sqrt{r^2 - \epsilon} = \frac{r}{\sqrt{r^2 - \epsilon}} > 1
\] 
by $\Phi_\epsilon$. The density changes from $r^{k}$ to $(r^2 - \epsilon)^{k/2}$, 
scaled by $\lambda^{-k}$. So the weighted length is scaled by $\lambda^{1-k}$, 
which is less than 1 because by hypothesis $k>1$. 
\end{proof}

Now we use the map $\Phi_\epsilon$ to show that the cluster is connected for certain 
densities.

\begin{proposition}
In the plane with density $r^p$, $p < -2$, any isoperimetric cluster (including 
the interior) must be connected and unbounded. 
\end{proposition}

\begin{proof}
We work in the Euclidean cone of \cref{prop:perimdensitymap} with only perimeter
density and origin corresponding to infinity back in the plane. 
For small enough $\epsilon$ we can apply the map $\Phi_\epsilon$ of 
\cref{lemma:areapreserved} to a component that does not contain the origin, 
yielding a cluster with the same area and less perimeter by 
\cref{lemma:lessperim}. Thus in the cone the cluster must be connected and 
contain the origin, which implies that back in the plane the cluster must be 
connected and unbounded. 
\end{proof}

\begin{remark}
Consider the plane with density $r^p$ for $p > 0$. Note that under the conformal map given in \cref{prop:perimdensitymap}, the perimeter density is always in 
the form $|w|^k$ for some $0 < k < 1$. Therefore this map does not 
decrease perimeter.  

Consider maps in the form $(r^k - \epsilon)^{1/k}$. The determinant of the 
Jacobian is given by 
\[
\det(D\varphi) = r^{k - 2} (r^k - \epsilon)^{2/k - 1}
\]
For $0 < k < 2$ and any $R > 0$, there exists a small enough $\epsilon$ such 
that $[\det(D\varphi)](r) < 1$ for all $r > R$. 
For $k > 2$ and any $R > 0$, there exists a small enough $\epsilon$ such that 
$[\det(D\varphi)](r) > 1$ for all $r > 1$. 
Therefore there are no area-increasing maps in this form that decrease
perimeter for perimeter density $r^p$ for $0 < p < 1$. 
\end{remark}

\section{Comparisons with Other Candidates}

In this section we focus on the plane and compare our standard double bubble with vertex at the origin 
against three other candidates, offering numerical and theoretical evidence that
our standard double bubble is isoperimetric. All candidates are in equilibrium 
and separate and enclose two regions of equal area $1$. Without loss of 
generality, they are plotted symmetric about the $y$-axis and shown here for 
density $r^2$.

\begin{figure}
    \centering
   \includegraphics[height=1.5in]{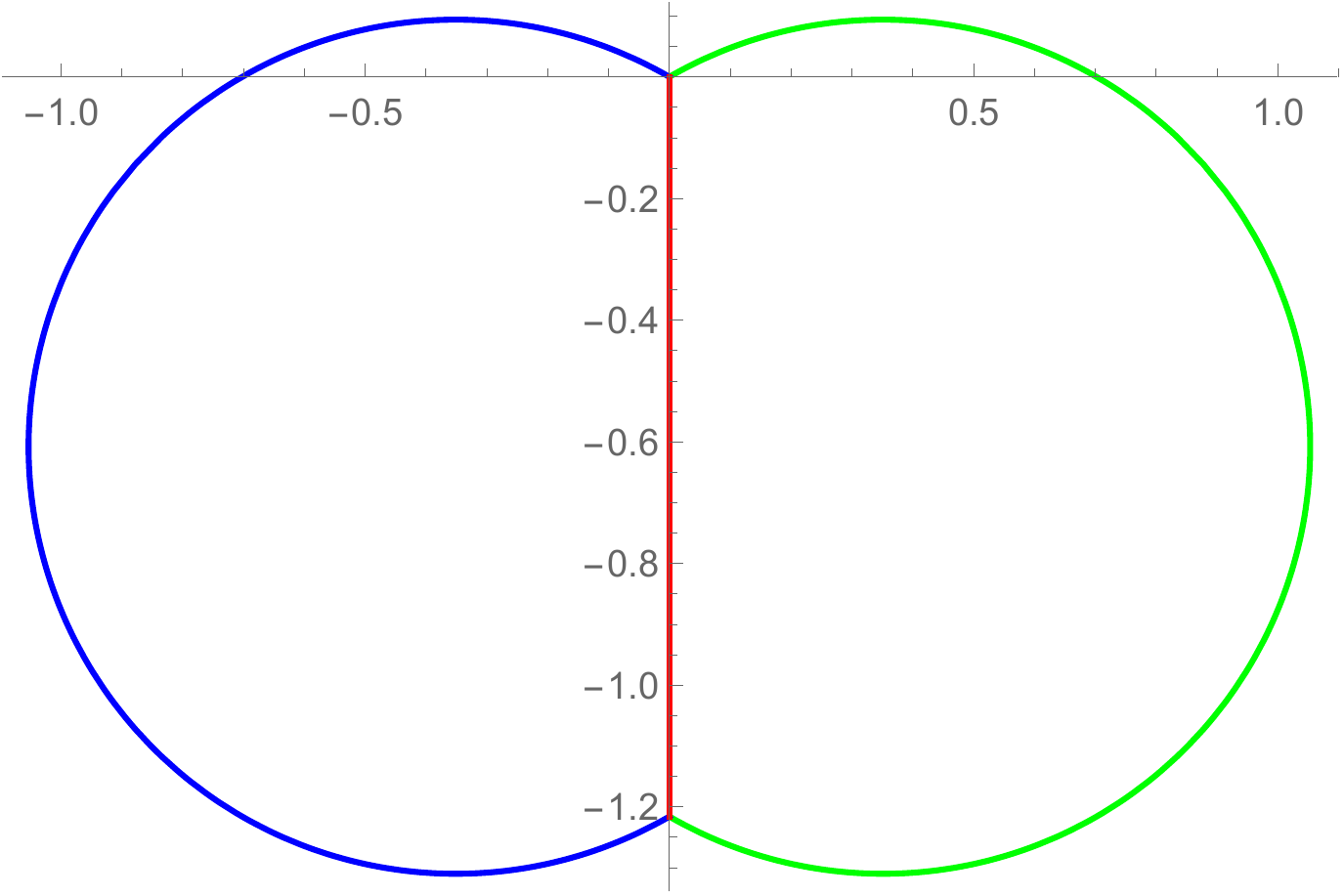}
    \caption{The standard double bubble, our conjectured champion, $p=2$.}
    \label{standard-double-bubble}
\end{figure}

\begin{figure}
    \centering
   \includegraphics[height=1.5in]{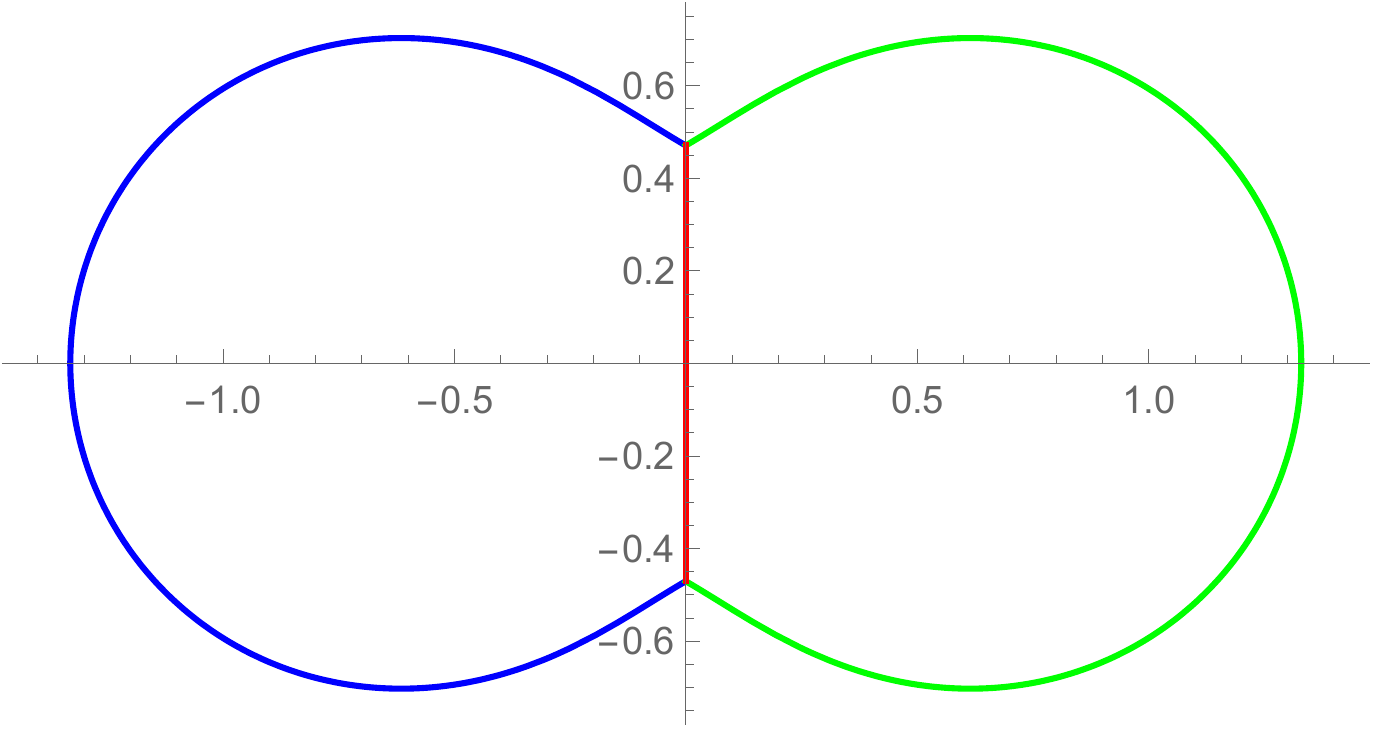}
    \caption{The symmetric double bubble, $p=2$.}
    \label{fig:symmetric}
\end{figure}

\begin{figure}
    \centering
   \includegraphics[height=1.5in]{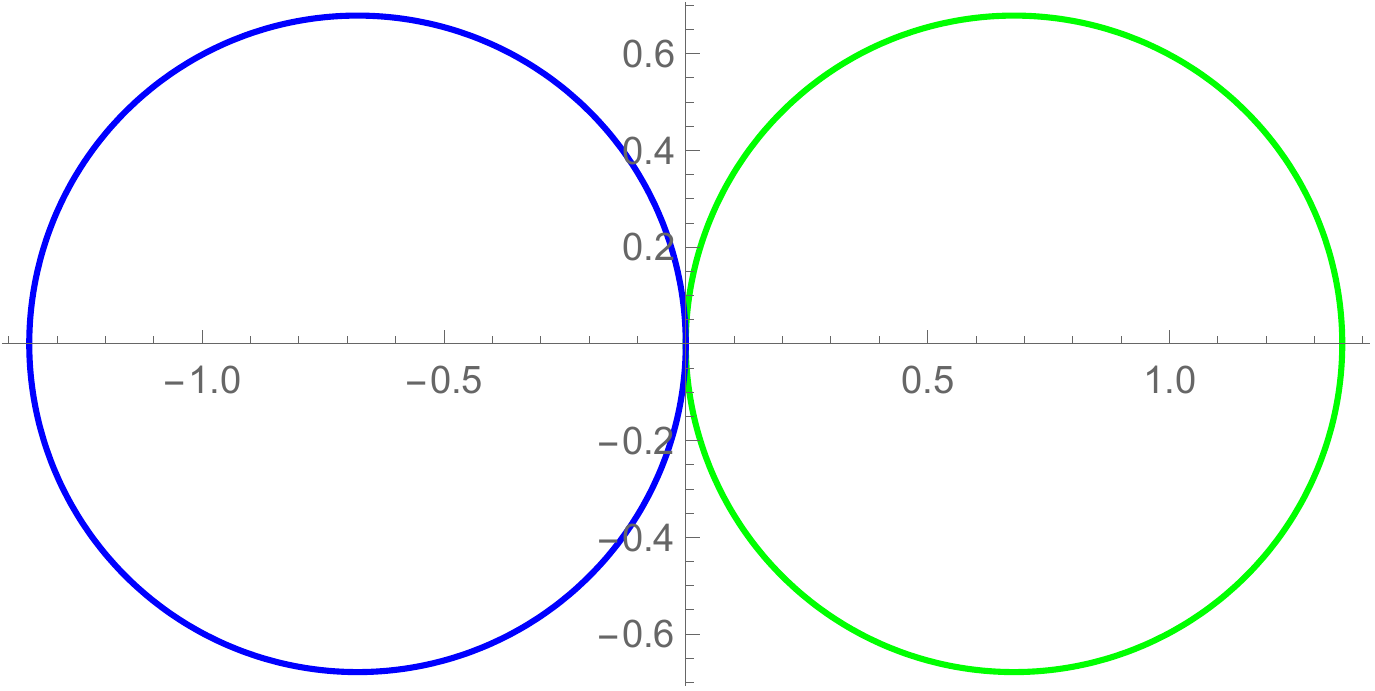}
    \caption{The two circles double bubble, $p=2$.}
    \label{fig:two-circles}
\end{figure}

\begin{figure}
    \centering
   \includegraphics[height=1.5in]{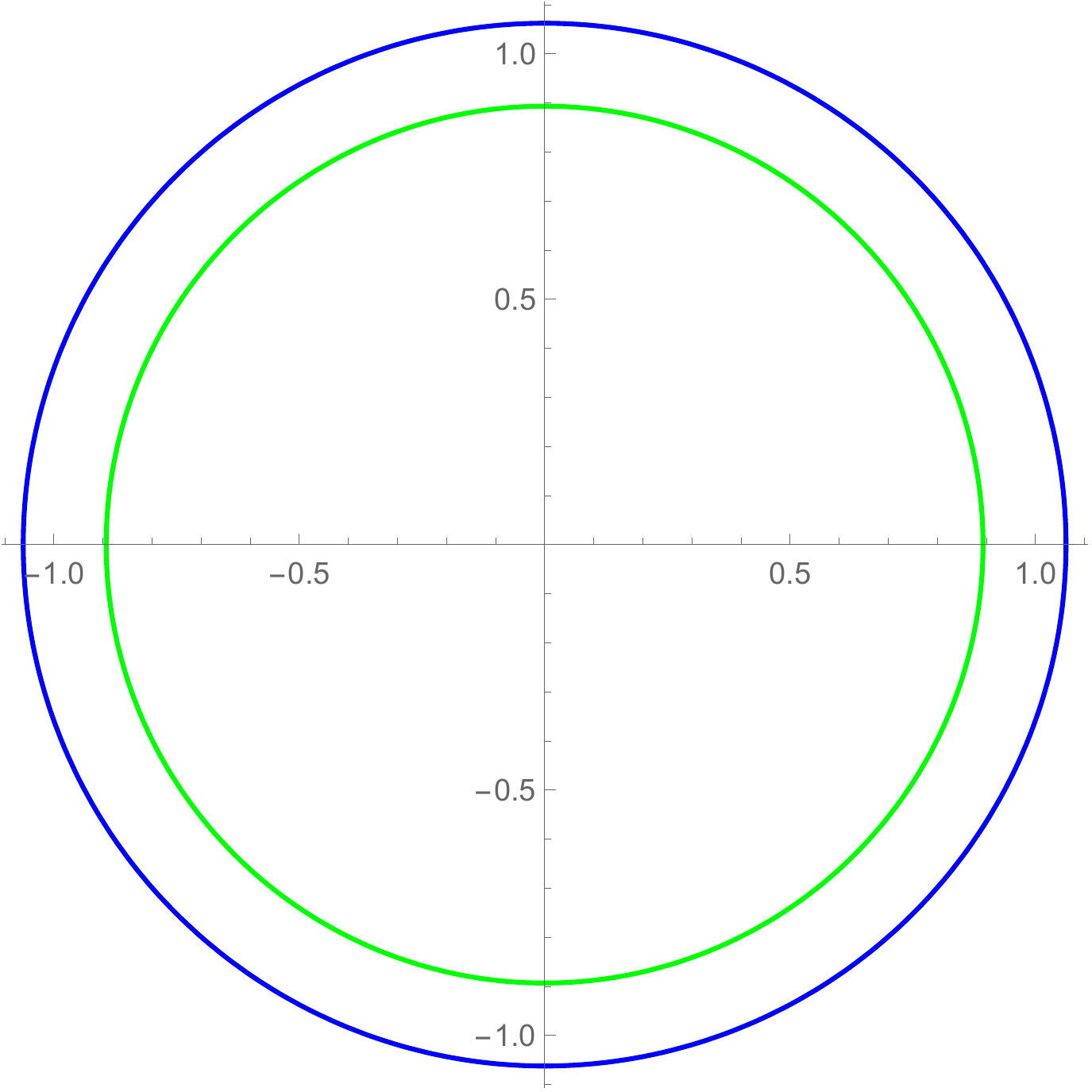}
    \caption{The concentric double bubble, $p=2$.}
    \label{fig:concentric}
\end{figure}

\cref{standard-double-bubble} shows our conjectured champion, the standard 
Euclidean double bubble with one vertex at the origin.
\cref{fig:symmetric} shows the next best candidate, a double bubble symmetric 
about the $y$-axis, composed of two constant-generalized-curvature arcs and a 
segment of the $y$-axis, meeting  at $120^\circ$. Note that the arcs do not have
constant \emph{Euclidean} curvature and hence are not circular.
\cref{fig:two-circles} shows the next best candidate, two circles meeting 
tangentially at the origin. Recall that a circle at the origin is the 
isoperimetric solution for the single bubble problem. Adding another circle, 
despite sharing no perimeter, does reasonably well, closely matching the 
perimeter of the symmetric double bubble for large $p$. The general $120^\circ$ 
equilibrium condition does not apply at the origin, because the density vanishes
there; indeed, \cref{sec:geodesics-in-plane} shows that shortest paths can have 
sharp (but not arbitrarily sharp) corners at the origin. Equilibrium still holds
for variations that are smooth diffeomorphisms, because each circle is 
minimizing. Nevertheless, \cref{prop:better-deformation} below shows that in 
fact equilibrium fails because perimeter can be reduced to first order by a 
Lipschitz deformation that pinches the two circles together, the very kind of 
deformation used in proving that curves meet at $120^\circ$ angles where the 
density is positive.
\cref{fig:concentric} shows two concentric circles, evidently worse even than 
the previous two circles candidate, because each of its bubbles does worse than 
a circle at the origin, the isoperimetric single bubble. Nonetheless, circles 
centered at the origin have constant generalized curvature, so the configuration
is in (unstable) equilibrium.
\cref{tab:perimeter-computation} gives the perimeters of the computed 
configurations in the plane with densities $r^p$, for $1 \le p \le 10.$

\begin{table}
    \centering\tiny\sffamily
    \begin{tabularx}{\textwidth}{X|cccccccccc}
        \toprule
        $p$ & $1$ & $2$ & $3$ & $4$ & $5$ & $6$ & $7$ & $8$ & $9$ & $10$ \\[1ex]
        standard & $6.490$ & $7.597$ & $8.979$ & $10.493$ & $12.085$ & $13.731$ 
        & $15.416$ & $17.132$ & $18.872$ & $20.632$ \\
        symmetric & $6.720$ & $7.837$ & $9.176$ & $10.650$ & $12.212$ & $13.835$
        & $15.502$ & $17.203$ & $18.932$ & $20.683$\\
        two circles & $6.868$ & $7.858$ & $9.177$ & $10.650$ & $12.212$ & 
        $13.835$ & $15.502$ & $17.203$ & $18.932$ & $20.683$\\
        concentric & $9.931$ & $12.009$ & $14.346$ & $16.820$ & $19.379$ & 
        $21.998$ & $24.661$ & $27.359$ & $30.085$ & $32.834$  \\
        \bottomrule
    \end{tabularx}
    \bigskip
    \caption{Perimeters of equilibrium double bubble candidates, rounded to the 
             nearest thousandth. Computations are done numerically in 
             Mathematica.}
    \label{tab:perimeter-computation}
\end{table}

\smallskip
The following proposition shows that although the candidate of 
\cref{fig:two-circles} is in equilibrium under smooth diffeomorphisms (because 
each circle is minimizing), it is not in equilibrium under small Lipschitz 
deformations about the origin that can pinch pieces together.

\begin{proposition}\label{prop:better-deformation}
A double bubble consisting of two circles tangent to each other at the origin is
not in equilibrium under (small) area-preserving Lipschitz deformations. 
\end{proposition}

\begin{figure}\label{tangentcirc}
    \centering
    \includegraphics{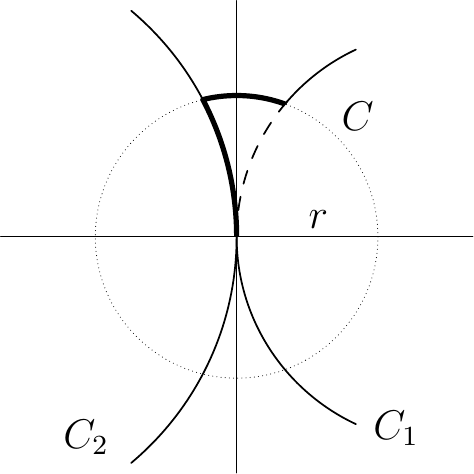}
    \caption{Deforming the dashed portion of $C_1$ onto the bold portion of $C$ 
             and $C_2$ reduces perimeter, belying equilibrium.}
    \label{fig:bubdotsl}
\end{figure}

\begin{proof}
Given small $\epsilon > 0$, there is a $\delta > 0$ such that part of the 
portion of the smaller circle $C_1$ in the lower half of an $\epsilon$-ball 
about the origin can be Lipschitz deformed to a chord to reduce area by any 
amount less than $\delta$.  As in \cref{fig:bubdotsl}, for small $r > 0$ first 
Lipschitz deform the top half of the arc of the smaller circle $C_1$ inside the 
circle $C$ about the origin of radius $r$ onto the other circle and a portion of
$C$. The perimeter saved is greater than the weighted length of a ray, which is 
$\int r^p \sim r^{p+1}$. The arc of $C$ adds perimeter on the order of 
$r^2 r^p$ = $r^{p+2}$. Hence for small $r$, perimeter is reduced. The area added
to the first region can now be returned by deforming an initial arc of $C_1$ 
below the origin to a chord, with further reduction of perimeter, all inside the
$\epsilon$-ball.
\end{proof}

\newpage
\printbibliography

@article{Morgan2003RegularityOI,
  title={Regularity of isoperimetric hypersurfaces in Riemannian manifolds},
  author={F. Morgan},
  journal={Transactions of the American Mathematical Society},
  year={2003},
  volume={355},
  pages={5041-5052}
}

@article{rosales,
author = {César Rosales and Antonio Cañete and Vincent Bayle and Frank 
          Morgan},
title = {On the isoperimetric problem in Euclidean space with density.},
fjournal = {Calculus of Variations and Partial
            Differential Equations},
journal = {Calc.\ Var.\ Partial Differential Equations},
year = {2008},
month = {1},
day = {01},
volume = {31},
number = {1},
pages = {27--46},
issn = {1432-0835},
doi = {10.1007/s00526-007-0104-y}
}

@article{foisy1993,
author = {Joel Foisy and Manuel Alfaro and Jeffrey Brock and Nickelous 
          Hodges and Jason Zimba},
fjournal = {Pacific Journal of Mathematics},
journal = {Pacific J.\ Math.},
issn = {0030-8730},
number = {1},
pages = {47--59},
publisher = {Pacific Journal of Mathematics, A Non-profit Corporation},
title = {The standard double soap bubble in {$\R^2$} uniquely 
         minimizes perimeter.},
url = {https://projecteuclid.org:443/euclid.pjm/1102634378},
volume = {159},
year = {1993}
}

@article{dahlberg2010,
author = {Jonathan Dahlberg and Alexander Dubbs and Edward Newkirk and Hung 
          Tran},
fjournal = {New York Journal of Mathematics},
journal = {New York J.\ Math.},
issn = {1076-9803},
pages = {31--51},
title = {Isoperimetric regions in the plane with density {$r^p$}.},
year = {2010},
volume = {16},
month = {5},
day = {05}
}

@article{G14,
author = {Wyatt Boyer and Bryan Brown and Gregory {R.\ Chambers} and Alyssa 
          Loving and Sarah Tammen},
title = {Isoperimetric regions in {$R^n$} with density {$r^p$}.},
year = {2016},
fjournal = {Analysis and Geometry in Metric Spaces},
journal = {Anal.\ Geom.\ Metr.\ Spaces},
pages = {236--265},
volume = {4},
issue = {1},
issn = {2299-3274},
eprinttype = {arxiv},
eprint = {1504.01720},
}

@unpublished{china19,
author = {Juiyu Huang and Xinkai Qian and Yiheng Pan and Mulei Xu and Lu 
          Yang and Junfei Zhou},
title = {Isoperimetric problems on the line with density {$\abs{x}^p$}.},
year = {2019},
month = {2},
day = {18},
url = {https://sites.williams.edu/Morgan/files/2019/02/x18Feb19.pdf},
addendum = {See \url{https://sites.williams.edu/Morgan/2019/02/22/china/}}
}

@book{morgan,
author = {Frank Morgan},
title = {Geometric Measure Theory: A Beginner's Guide},
month = {5},
year = {2016},
edition = {5},
publisher ={Academic Press},
isbn = {9780128044896}
}

@article{annals,
ISSN = {0003486X},
URL = {http://www.jstor.org/stable/3062123},
author = {Michael Hutchings and Frank Morgan and Manuel Ritoré and Antonio
          Ros},
fjournal = {Annals of Mathematics},
journal = {Ann.\ of Math.},
number = {2},
pages = {459--489},
publisher = {Annals of Mathematics},
title = {Proof of the double bubble conjecture.},
volume = {155},
year = {2002}
}

@article{hutchings,
author = {Michael Hutchings},
title = {The structure of area-minimizing double bubbles.},
fjournal = {The Journal of Geometric Analysis},
journal = {J.\ Geom.\ Anal.},
year = {1997},
month = {06},
day = {01},
volume = {7},
number = {2},
pages = {285--304},
issn = {1559-002X},
doi = {10.1007/BF02921724},
}

@article{G06,
title = {The isoperimetric problem on planes with density},
author = {Colin Carroll and Adam Jacob and Conor Quinn and Robin Walters},
journal = {Bull. Austral. Math. Soc.},
year = {2008},
volume = {78},
pages = {177--197},
doi = {10.1017/S000497270800052X}
}

@article{wichiramala,
author = {Wacharin Wichiramala},
title = {Proof of the planar triple bubble conjecture},
month = {02},
day = {05},
volume = {2004},
issue = {567},
pages = {1--49},
issn = {1435-5345},
journal = {J.\ Reine Angew.\ Math.},
doi = {10.1515/crll.2004.011}
}

@online{nardulli,
author = {Stefano Nardulli and Pierre Pansu},
title = {A discontinuous isoperimetric profile for a complete Riemannian 
         manifold},
year = {2015},
month = {06},
day = {17},
archivePrefix = {arXiv},
eprint = {1506.04892}
}

@article{papasoglu,
author = {Panos Papasoglu and Eric Swenson},
title = {A surface with discontinuous isoperimetric profile and expander 
         manifolds},
archivePrefix = {arXiv},
eprint = {1509.02307},
year = {2019},
month = {07},
day = {21},

}

@article{Morgan2013,
author="Morgan, Frank
and Pratelli, Aldo",
title="Existence of isoperimetric regions in $\R^n$ with density",
journal="Annals of Global Analysis and Geometry",
year="2013",
month="4",
day="01",
volume="43",
number="4",
pages="331--365",
issn="1572-9060",
doi="10.1007/s10455-012-9348-7",
url="https://doi.org/10.1007/s10455-012-9348-7"
}

@article{Reichardt2007,
author="Reichardt, Ben W.",
title="Proof of the Double Bubble Conjecture in $\R^n$",
journal="Journal of Geometric Analysis",
year="2007",
month="11",
day="16",
volume="18",
number="1",
pages="172",
issn="1559-002X",
doi="10.1007/s12220-007-9002-y",
url="https://doi.org/10.1007/s12220-007-9002-y"
}

\end{document}